\documentclass[a4paper,10pt,reqno]{amsart}
\usepackage{amsfonts,amsmath,amssymb,bbm,amsthm,amsbsy}
\usepackage{graphicx}
\usepackage{geometry}
\usepackage{float}
\usepackage{verbatim}
\usepackage{setspace}
\usepackage{siunitx}
\usepackage{paralist}
\usepackage{enumitem}
\usepackage{xcolor}
\usepackage{lipsum}
\usepackage[english]{babel}
\usepackage{booktabs}
\usepackage{array}
\usepackage{subfig}
\usepackage{caption}
\usepackage{mathtools}
\usepackage{stmaryrd}
\usepackage{stackengine,scalerel}
\usepackage{hyperref}
\usepackage[utf8]{inputenc}
\usepackage{relsize}

\newenvironment{frcseries}{\fontfamily{frc}\selectfont}{}

\makeatletter
\newtheorem*{rep@theorem}{\rep@title}
\newcommand{\newreptheorem}[2]{
\newenvironment{rep#1}[1]{
 \def\rep@title{#2 \ref{##1}}
 \begin{rep@theorem}}
 {\end{rep@theorem}}}
\makeatother

\newtheorem{thm}{Theorem}[section]
\newtheorem{lemma}[thm]{Lemma}
\newtheorem{prop}[thm]{Proposition}

\newtheorem*{thm*}{Theorem}
\newtheorem*{lemma*}{Lemma}
\newtheorem*{prop*}{Proposition}
\newtheorem*{corr*}{Corrolary}
\newtheorem*{claim*}{Claim}

\theoremstyle{remark}
\newtheorem{rmk}[thm]{Remark}

\newtheorem*{rmk*}{Remark}
\newtheorem*{conj*}{Conjecture}
\newtheorem*{quest*}{Question}

\theoremstyle{definition}
\newtheorem{defn}[thm]{Definition}
\newtheorem{exmp}[thm]{Example}

\newtheorem*{defn*}{Definition}
\newtheorem*{exmp*}{Example}

\newreptheorem{theorem}{Theorem}
\newreptheorem{corollary}{Corollary}
\newreptheorem{proposition}{Proposition}

\usepackage{fancyhdr}
\usepackage{tikz}
\usepackage{tikz-cd}

\usepackage{IEEEtrantools}

\newenvironment{equ*}[1]{\begin{IEEEeqnarray*}{#1}}{\end{IEEEeqnarray*}}

\newcommand{\crotimes}{\,\,{\widehat{\otimes}}_\mathit{R}\,\,}

\newcommand{\Z}{\mathbb{Z}}

\DeclareFontFamily{U}{mathx}{}
\DeclareFontShape{U}{mathx}{m}{n}{<-> mathx10}{}
\DeclareSymbolFont{mathx}{U}{mathx}{m}{n}
\DeclareMathAccent{\widecheck}{0}{mathx}{"71}

\newcommand{\inj}{\hookrightarrow}
\newcommand{\sur}{\twoheadrightarrow}

\newcommand{\kG}{k\llbracket G\rrbracket}
\newcommand{\RG}{R\llbracket G\rrbracket}

\newcommand{\kH}{k\llbracket H\rrbracket}

\DeclareMathOperator{\Hom}{Hom}
\DeclareMathOperator{\Cod}{Cod}
\DeclareMathOperator{\Fun}{Fun}
\DeclareMathOperator{\IFun}{IFun}

\DeclareMathOperator{\Tor}{Tor}

\DeclareMathOperator{\Ind}{Ind}

\DeclareMathOperator{\Res}{Res}

\newcommand{\Set}{\mathbf{Set}}
\newcommand{\Grp}{\mathbf{Grp}}
\newcommand{\Ab}{\mathbf{Ab}}

\newcommand{\ModR}{\mathbf{Mod}(R)}

\newcommand{\Pro}{\mathbf{Pro}}

\newcommand{\BPro}{\mathbf{BPro}}
\newcommand{\BProf}{\mathbf{BSet}_{\mathrm{fin}}}
\newcommand{\Cat}{\mathbf{Cat}}
\newcommand{\Graph}{\mathbf{Graph}}
\newcommand{\PGraph}{\mathbf{PGraph}}

\newcommand{\PGrp}{\mathbf{PGrp}}
\newcommand{\BPGrp}{\mathbf{BPGrp}}
\newcommand{\BPGrpf}{\mathbf{BGrp}_{\mathrm{fin}}}

\newcommand{\PAb}{\mathbf{PAb}}

\newcommand{\PModR}{\mathbf{PMod}(R)}
\newcommand{\BPModR}{\mathbf{BPMod}(R)}
\newcommand{\BPModRf}{\mathbf{BMod}(R)_{\mathrm{fin}}}

\newcommand{\coshR}{\mathbf{CoSh}(R)}

\newcommand{\G}{\mathcal{G}}

\title{Coproducts Internal to Profinite Spaces}
\author{Jiacheng Tang}
\thanks{Email: jiacheng.tang@postgrad.manchester.ac.uk; Address: Alan Turing Building, University of Manchester, Manchester, M13 9PY, United Kingdom.}

\calclayout
\begin{document}
\maketitle

\begin{abstract}
We give a categorical explanation for many properties of profinite coproducts of profinite groups, which were previously proven on a case-by-case basis. All of these properties take the form ``certain functors preserve profinite coproducts". We give various examples to show how our framework can be applied. We also point out connections to internal category theory and profinite Bass--Serre theory.
\end{abstract}

\section{Introduction}
\label{sec1}

The category $\PGrp$ of profinite groups has all small coproducts, but infinite coproducts are difficult to use in practice. On the other hand, \cite{haran}, \cite{melfree} and \cite{ribesgraph} develop a notion of infinite coproducts of profinite groups which are indexed not by sets, but by profinite spaces. This has numerous applications, as illustrated in those papers. There has been recent work, such as \cite{gareth}, \cite{boggi} and \cite{jc3}, which presents more categorical approaches to studying this notion of topological coproducts. However, there are still many results on these topological coproducts in the literature which are proven on a case-by-case basis and do not yet seem to have any satisfactory categorical explanation. To demonstrate our point, consider the following. Let $\coprod_{x\in X}G_x$ be the ``topological coproduct" of profinite groups $G_x$ indexed by a profinite space $X$, and let $\widehat{\bigoplus}_{x\in X}M_x$ be the corresponding notion for profinite modules $M_x$ over a fixed profinite ring $R$. These will be defined properly later (see Example \ref{egdisc}). We now state a few facts about these coproducts and put what the right categorical explanation should be in brackets at the end.

\begin{itemize}
\item (Example \ref{maineg}\ref{maineg4}) If $G^\mathrm{ab}$ denotes the abelianisation of a profinite group $G$, then $$\left(\coprod_{x\in X}G_x\right)^\mathrm{ab}=\widehat{\bigoplus}_{x\in X}(G_x^\mathrm{ab}).$$ In particular, if $I$ is a set, then the abelianisation of the free profinite group on $I$ (converging to 1) is the free proabelian group on $I$ (converging to 1).
\\(Abelianisation commutes with coproducts.)

\item (\cite[Theorem 9.1.1]{ribesgraph}) If $N$ is a profinite left $R$-module, then $$\Tor^R_i\left(\widehat{\bigoplus}_{x\in X}M_x,N\right)=\widehat{\bigoplus}_{x\in X} \Tor^R_i(M_x,N).$$
($\Tor$ commutes with coproducts.)

\item (\cite[Proposition A.14]{garethrel}) If $f\colon Y\sur X$ is a surjection of profinite spaces, then $$R\llbracket Y\rrbracket=\widehat{\bigoplus}_{x\in X}R\llbracket f^{-1}(x)\rrbracket.$$
(The free profinite $R$-modules functor $R\llbracket-\rrbracket$ commutes with coproducts.)
\end{itemize}

\subsection*{Main Aim} The aim of this paper is to provide, in Theorem \ref{mainthm}, a general categorical framework that explains all of the above and more. In fact, \cite[Proposition A.13]{garethrel} already suggests that categorical explanations should exist for the above facts, and our main theorem can be viewed as a generalisation of \cite[Proposition A.13]{garethrel}. Informally, we will show that many naturally occurring functors in the profinite setting (many of which are automatically left adjoints) commute with these topological coproducts. These coproducts have various existing applications already: \cite[Theorem 9.4.1]{ribesgraph} uses them to construct a Mayer--Vietoris sequence for profinite groups acting on trees, \cite[Theorem 2.19]{jc3} uses them to prove Mackey's Formula for general profinite groups, and so on. All of the stated results mimic classical results involving ordinary groups and modules, but use this notion of topological coproducts in place of ordinary coproducts. The generality of our main theorem (Theorem \ref{mainthm}), which we illustrate in Example \ref{maineg}, thus provides a concrete strategy for transporting classical results of categorical flavours into the profinite world.

\subsection*{Outline of Paper}

The main body of this paper is divided into Sections \ref{sec2} and \ref{sec3}, which are somewhat disjoint. Readers who are only interested in the main theorem mentioned above could therefore skip Section \ref{sec2} altogether. In Section \ref{2.1}, we recall the definition of an internal category (specifically internal to $\Pro$, the category of profinite spaces), as it turns out that the topological coproducts mentioned above are, in a precise sense, a kind of internal coproducts (Example \ref{egdisc}). More generally, we can define colimits internal to $\Pro$, which are constructed from coproducts and coequalisers (Theorem \ref{thm1}). We will see, in Section \ref{2.2}, that internal pushouts of profinite groups in our sense are the profinite fundamental groups of certain graphs of groups, in exactly the same way as classical Bass--Serre theory.

In Section \ref{sec3}, we focus on internal coproducts (rather than all colimits) and study when they are preserved by functors. The conditions we need are mild (Proposition \ref{commprop}) and are often automatic, since many naturally occurring functors in the profinite setting arise from \emph{relative adjunctions} (Lemma \ref{lem1} and Example \ref{exampleJco}). In this case, we get more than just the preservation of coproducts (Proposition \ref{leftadjgood}). We will end the paper with Example \ref{maineg}, which shows that there is an abundance of situations where our theorem applies. In the appendix, we will briefly discuss Pontryagin duality in our context.

We will assume basic knowledge of category theory (refer to \cite{maccat}) and profinite groups (refer to \cite{profinite}). Section \ref{sec2} will require knowledge of internal category theory (refer to \cite{handbook}) and profinite graphs (refer to \cite{ribesgraph}), while Section \ref{sec3} will require knowledge of pro-completions (refer to \cite{elephant}).

Remark: When we write ``=" in this paper, we usually mean canonically isomorphic (or equivalent in the case of categories).

Convention: Unless otherwise specified, all rings are associative with a 1 but not necessarily commutative.

\subsection*{Acknowledgements}
The author would like to thank his supervisor Peter Symonds for his constant guidance, Calum Hughes and Giacomo Tendas for helpful discussions, as well as Giacomo Tendas for reading drafts of this paper and giving useful feedback. The author would also like to thank the anonymous referee for their useful comments.

This work was supported by the EPSRC Doctoral Training Partnership [grant number \linebreak EP/W524347/1].

\section{Categories Internal to $\Pro$}
\label{sec2}

In this section, we will recall the concept of internal categories, following mostly the notations of \cite[Chapter 8]{handbook}. We will see that the topological coproducts of profinite groups appearing in the introduction can be described as coproducts internal to $\Pro$ (Example \ref{egdisc}). More generally, colimits internal to $\Pro$ can be constructed from coproducts and coequalisers (Theorem \ref{thm1}), as one would expect. In Section \ref{2.2}, we will show that our notion of colimits internal to $\Pro$ in fact occurs naturally in the context of profinite Bass--Serre theory.

\subsection{Internal Categories and Colimits}\label{2.1} We recall the following special case of \cite[Chapter 8]{handbook}. Let $\Pro$ denote the category of profinite spaces (and continuous maps) and let $A$ be a category internal to $\Pro$. Spelt out, this means that there are two spaces $A_0, A_1\in\Pro$ and maps $$d_0, d_1\colon A_1\to A_0,\,\,\,\, i\colon A_0\to A_1,\,\,\,\, c\colon A_1\times_{A_0}A_1\to A_1$$ in $\Pro$, where the pullback is formed using $(d_1,d_0)$, such that the underlying sets and functions describe an ordinary category. One can then define functors and natural transformations internal to $\Pro$. Let $\Cat(\Pro)$ denote the ordinary 1-category of categories and functors internal to $\Pro$.

Even though $\Pro$ is not a category internal to itself, it makes sense to talk about internal functors from an internal category $A$ to $\Pro$. Spelt out, an internal functor $P\colon A\to\Pro$ consists of a space $P_0\in\Pro$ and maps $$p_0\colon P_0\to A_0,\,\,\,\, p_1\colon A_1\times_{A_0}P_0\to P_0,$$ where the pullback is formed using $(d_0,p_0)$, such that the underlying sets and functions describe an ordinary functor $A\to\Set$ which sends an object $a\in A_0$ to the set $P(a)=p_0^{-1}(\{a\})$, and which sends a map $f\colon a\to b$ in $A_1$ to the function of sets $$P(f)\colon P(a)\to P(b),\,\,\,\, x\mapsto p_1(f,x).$$ Given two such functors $P, Q\colon A\to\Pro$, an internal natural transformation is a map $\alpha\colon P_0\to Q_0$ in $\Pro$, such that the underlying sets and functions describe an ordinary natural transformation. Internal functors from $A$ to $\Pro$ and their natural transformations form an ordinary category, which we denote by $\IFun(A,\Pro)$. The reader should consult \cite{handbook} for the omitted details.

There is an ordinary functor $\Delta_A\colon \Pro\to\IFun(A,\Pro)$ which sends $X\in\Pro$ to the constant internal functor $\Delta_A(X)_0=X\times A_0\to A_0$. As $\Pro$ has finite limits and coequalisers, it is easy to show that $\Delta_A$ has a left adjoint for every $A\in\Cat(\Pro)$, i.e.\ that $\Pro$ has all internal colimits (see \cite[Proposition 8.3.2]{handbook}). We point out that $\Pro$ does not have all internal limits, since it is not Cartesian closed (see \cite[Proposition 8.3.5]{handbook}).

To motivate what follows, let us forget that we are dealing with internal categories for a moment and think about $\Set$. Given any $A\in\Cat$, the functor $\Delta_A\colon\Set\to\Fun(A,\Set)$ has a left adjoint, i.e.\ $\Set$ has all (small) colimits. Since both $\Set$ and $\Fun(A,\Set)$ have finite products, we can look at group objects (or ring objects, or $R$-module objects...) in them. The functor $\Delta_A$ preserves finite products, so induces a functor on group objects (also denoted by $\Delta_A$), which is simply the functor $\Delta_A\colon\Grp\to\Fun(A,\Grp)$ that sends $G\in\Grp$ to the constant functor with value $G$. This new functor $\Delta_A$ on $\Grp$ also has a left adjoint for any $A\in\Cat$, i.e.\ $\Grp$ has all (small) colimits. Note that this does not follow automatically from the fact that $\Set$ has all colimits and is genuinely a theorem about $\Grp$.

Now, let us return to categories internal to $\Pro$. It is easy to check that given $A\in\Cat(\Pro)$, the internal functor category $\IFun(A,\Pro)$ has all finite products. The terminal object is given by the space $A_0$. Given two internal functors $P,Q\in\IFun(A,\Pro)$, their product $P\times Q$ is the internal functor with space $(P\times Q)_0=P_0\times_{A_0}Q_0$, and where $(p\times q)_0\colon P_0\times_{A_0}Q_0\to A_0$ is the obvious map, and $(p\times q)_1\colon A_1\times_{A_0}(P_0\times_{A_0}Q_0)\to P_0\times_{A_0}Q_0$ is the unique map induced by
\[\begin{tikzcd}
A_1\times_{A_0}P_0\times_{A_0}Q_0 \arrow[rdd, "p_1\circ\pi"', bend right] \arrow[rrd, "q_1\circ\pi", bend left] \arrow[rd, dashed] &                                                       &                       \\
                                                                                                                                   & P_0\times_{A_0}Q_0 \arrow[d, "\pi"] \arrow[r, "\pi"'] & Q_0 \arrow[d, "q_0"'] \\
                                                                                                                                   & P_0 \arrow[r, "p_0"]                                  & A_0                  
\end{tikzcd}\]
(Note that we are using the same letter $\pi$ to denote different projection maps, but there should be no confusion.) The careful reader might wish to check that it is the binary product in $\IFun(A,\Pro)$, which is not difficult.

Similar to our previous discussion for sets and groups, let us now consider group objects in $\Pro$ and $\IFun(A,\Pro)$. The functor $\Delta_A\colon\Pro\to\IFun(A,\Pro)$ again preserves finite products, so induces a functor $\Delta_A\colon\PGrp\to\Grp(\IFun(A,\Pro))$, where $\PGrp=\Grp(\Pro)$ is the category of profinite groups. It sends $G\in\PGrp$ to the constant functor $(\Delta_A(G))_0=G\times A_0\to A_0$. If $P\in\Grp(\IFun(A,\Pro))$ has a reflection along $\Delta_A$, then we call the reflection the \emph{colimit of $P$ (internal to $\Pro$)}. If such a colimit exists for a fixed $A$ and for every $P$, i.e.\ if $\Delta_A$ has a left adjoint, then we say the category $\PGrp$ \emph{has all colimits of shape $A$ internal to $\Pro$}. If such a left adjoint exists for every $A$, then we say that the category $\PGrp$ \emph{has all colimits internal to $\Pro$}. We will soon show that this is in fact the case (see Theorem \ref{thm1}). We have chosen to work specifically with group objects here for concreteness, but our main theorem (Theorem \ref{mainthm}) will be stated more generally.

Before we look at some examples of these colimits, let us make explicit what the objects and maps in $\Grp(\IFun(A,\Pro))$ look like. An object $P\in\Grp(\IFun(A,\Pro))$ consists of:
\begin{itemize}
\item a profinite space $P_0$, which as a set is a disjoint union $\bigsqcup_{a\in A_0}P(a)$;
\item a map $p_0\colon P_0\to A_0$ in $\Pro$ which sends $P(a)$ to $a$;
\item a map $p_1\colon A_1\times_{A_0} P_0\to P_0$ in $\Pro$ which sends $(f,x)\in A_1\times_{A_0}P_0$ (where $d_0(f)=a$ and $x\in P(a)$) to an element $P(f)(x)\in P(b)$ (where $d_1(f)=b$);
\item maps $P_0\times_{A_0}P_0\to P_0$, $P_0\to P_0$ and $A_0\to P_0$ in the slice category $\Pro_{/A_0}$ (which we think of as multiplication, inversion and identity);
\item such that when restricted to each $P(a)=p_0^{-1}(\{a\})$, these three maps make $P(a)$ a group;
\item and such that the assignment $a\mapsto P(a)$, $f\mapsto P(f)$ describes a functor $A\to\Grp$.
\end{itemize}
Given two such objects $P$ and $Q$, a map $\alpha\colon P\to Q$ in $\Grp(\IFun(A,\Pro))$ consists of:
\begin{itemize}
\item a map $\alpha\colon P_0\to Q_0$ in $\Pro_{/A_0}$;
\item such that the restrictions $\alpha_a\colon P(a)\to Q(a)$ describe a natural transformation from $P$ to $Q$, viewed as functors $A\to\Grp$.
\end{itemize}

\begin{exmp}\label{egdisc}
\begin{enumerate}[label=(\roman*)]
\item\label{intcolim1} Let $A$ be an ordinary finite category, i.e.\ $A_0$ and $A_1$ are both finite. Then, $A$ is canonically a category internal to $\Pro$ by equipping $A_0$ and $A_1$ each with the discrete topology.

In this case, an internal functor $A\to\Pro$ is just an ordinary functor etc., so the functor $\Delta_A\colon\PGrp\to\Grp(\IFun(A,\Pro))=\Fun(A,\PGrp)$ does have a left adjoint, since the category $\PGrp$ has all (ordinary) finite colimits.
\item\label{intcolim2} Let $X\in\Pro$ be a profinite space. Then, we canonically obtain a category internal to $\Pro$, which we denote by $X^\delta$, called the \emph{discrete internal category associated to $X$}. It is given by $X^\delta_0=X$, $X^\delta_1=X$, and all the maps $d_0, d_1, i, c$ are identity maps. (Note that ``discrete" here is referring not to the topology of $X$, but to the fact that the only morphisms in $X^\delta$ are the identity maps on objects.)

In this case, the internal functor category $\IFun(X^\delta,\Pro)$ is equivalent to the slice category $\Pro_{/X}$ (see \cite[Proposition 8.2.5]{handbook}). An object $P\in\Grp(\IFun(X^\delta,\Pro))$ is exactly what \cite[Section 5.1]{ribesgraph} calls a \emph{sheaf of profinite groups over $X$}. The functor $$\Delta_{X^\delta}\colon\PGrp\to\Grp(\IFun(X^\delta,\Pro))$$ does have a left adjoint by \cite[Proposition 5.1.2]{ribesgraph}. Following \cite{ribesgraph}, we will call the image of $P\in\Grp(\IFun(X^\delta,\Pro))$ under the left adjoint the \emph{free (profinite) product of $P$} and denote it by $\coprod_{x\in X}P(x)\in\PGrp$, where the topology on $P_0=\bigsqcup_{x\in X} P(x)$ is implicit. Note that this is the ``topological coproduct" mentioned in the introduction (except we now write $P(x)$ instead of $P_x$). Notice that, while the free product of ordinary groups is the coproduct in $\Grp$, the free profinite product of profinite groups, as defined above, is really a ``coproduct internal to $\Pro$", rather than the categorical coproduct in $\PGrp$.
\end{enumerate}
\end{exmp}

\begin{rmk}\label{rmk2}
Let $R$ be a profinite ring. If we replace groups with $R$-modules in \ref{intcolim2} above, interpreted appropriately (see below), then we get the \emph{profinite direct sum} $\widehat{\oplus}$ of profinite modules, which is defined in \cite[Section 9.1]{ribesgraph} and is the main object of study in \cite{jc3}. This is the ``topological coproduct" of profinite modules mentioned in the introduction.

To elaborate, note that the functor $\Delta_A\colon\Pro\to\IFun(A,\Pro)$ induces a functor on ring objects. Given a profinite ring $R$, let $R\times A_0$ denote its image under $\Delta_A$. Now observe that $\Delta_A$ induces a functor from the category of $R$-module objects in $\Pro$ (i.e.\ profinite $R$-modules) to the category of $(R\times A_0)$-module objects in $\IFun(A,\Pro)$. This functor, still denoted by $\Delta_A$, sends a profinite $R$-module $M$ to the constant functor $M\times A_0$. If $X\in\Pro$, then the left adjoint to $\Delta_{X^\delta}$ is the profinite direct sum $\widehat{\oplus}_X$.
\end{rmk}

As we saw in \ref{intcolim1} and \ref{intcolim2} of the example above, the category $\PGrp$ has finite colimits (so in particular coequalisers) and all coproducts internal to $\Pro$. The obvious next question to ask is whether this implies $\PGrp$ has all colimits internal to $\Pro$, and the answer is yes, as we shall now show.

\begin{thm}\label{thm1}
Let $A\in\Cat(\Pro)$. If $P\in\Grp(\IFun(A,\Pro))$, then its colimit (internal to $\Pro$) exists and can be constructed using coequalisers and coproducts of shapes $A_0^\delta$ and $A_1^\delta$. In particular, $\PGrp$ has all colimits internal to $\Pro$.
\begin{proof}
This is similar to the proof of the classical statement (cf.\ \cite[Theorem V.2.1]{maccat}), but we have to work in the internal language of $\Pro$. In fact, we will see below (Remark \ref{rmk1}) that this is just a special case of \cite[Lemma 7.4.8]{catlogic}, but we include a proof anyway for completeness\footnote{In fact, cocompleteness.}.

Recall that $A_i^\delta$ is the discrete internal category and that $\IFun(A_i^\delta,\Pro)=\Pro_{/A_i}$. Define an object $Q\in\Grp(\Pro_{/A_0})$, where $Q_0=P_0$ and $q_0=p_0\colon Q_0\to A_0$, and another object $R\in\Grp(\Pro_{/A_1})$, where $R_0=A_1\times_{A_0}P_0$ (the pullback of $d_0$ and $p_0$) and $r_0=\pi_{A_1}\colon R_0\to A_1$ is the canonical projection. Note that each fibre $R(f)=r_0^{-1}(\{f\})=P(d_0(f))$ (where $f\in A_1$) is indeed a group. Let $\coprod_{a\in A_0}P(a)$ and $\coprod_{f\in A_1}P(d_0f)$ be the free products of $Q$ and $R$ respectively.

We shall now define two maps $\phi,\psi\colon\coprod_{f\in A_1}P(d_0f)\to\coprod_{a\in A_0}P(a)$ in $\PGrp$, which is the same as defining maps $R_0=A_1\times_{A_0}P_0\to\coprod_{a\in A_0}P(a)$ in $\Pro$ such that their restriction to each fibre is a group homomorphism $R(f)=P(d_0f)\to\coprod_{a\in A_0}P(a)$ (where $f\in A_1$). The map $\phi$ is given by $$\phi\colon R_0=A_1\times_{A_0}P_0\overset{\pi}{\longrightarrow} P_0=Q_0\overset{\eta_Q}{\longrightarrow}\left(\coprod_{a\in A_0}P(a)\right)\times A_0\overset{\pi}{\longrightarrow}\coprod_{a\in A_0}P(a),$$ where $\eta$ is the unit of the adjunction involving $\Delta_{A_0^\delta}$. The restriction of $\phi$ to the fibre above $f\in A_1$ is the canonical map $P(d_0f)\to\coprod_{a\in A_0}P(a)$. On the other hand, the map $\psi$ is given by $$\psi\colon R_0=A_1\times_{A_0}P_0\overset{p_1}{\longrightarrow} P_0=Q_0\overset{\eta_Q}{\longrightarrow}\left(\coprod_{a\in A_0}P(a)\right)\times A_0\overset{\pi}{\longrightarrow}\coprod_{a\in A_0}P(a).$$ The restriction of $\psi$ to the fibre above $f\in A_1$ is the map $P(d_0f)\overset{P(f)}{\longrightarrow} P(d_1f)\to\coprod_{a\in A_0}P(a)$.

Let $C\in\PGrp$ be the coequaliser of $\phi$ and $\psi$. We leave it to the reader to check that $C$ is the colimit of $P$ (internal to $\Pro$).
\end{proof}
\end{thm}

\begin{rmk}\label{rmk1}
The above theorem is simply a special case of \cite[Lemma 7.4.8]{catlogic}. The reader should refer to the cited book for the following discussion.

Consider the category $\BPGrp$ whose objects are pairs $(P,X)$ with $X\in\Pro$ and $P \in\Grp(\Pro_{/X})$\footnote{Readers who are familiar with the Grothendieck construction will recognise this as a special instance of the said construction.}. A morphism from $(P,X)$ to $(Q,Y)$ in $\BPGrp$ consists of continuous maps $f\colon P\to Q$ and $X\to Y$ making the obvious square commute, such that the restriction of $f$ to each fibre above $X$ is a group homomorphism (see \cite[Section 5.1]{ribesgraph}; the notation $\BPGrp$ stands for ``bundles of profinite groups".). There is a functor $\Cod\colon \BPGrp\to\Pro$, called the \emph{codomain functor}, which sends an object $(P,X)\in \BPGrp$ to its codomain $X$, and it is easy to see that $\Cod$ is a Cartesian fibration (cf.\ \cite[Proposition 1.1.6]{catlogic}).

As pointed out in \cite[Definition 7.4.1]{catlogic}, an object of $\Grp(\IFun(A,\Pro))$ in our sense is exactly what \cite{catlogic} would call an \emph{internal diagram of type $A$ in $\Cod$}. In fact, $\Grp(\IFun(A,\Pro))$ is exactly the full subcategory of the category $\BPGrp^A$ in \cite[Definition 7.4.5]{catlogic} consisting of objects with $I=*$. As such, the above theorem is just (the dual of) \cite[Lemma 7.4.8]{catlogic} and its proof restricted to $I=*$.
\end{rmk}

\subsection{Connections to Profinite Bass--Serre Theory}\label{2.2} Free profinite products and finite colimits of profinite groups are used extensively in profinite Bass--Serre theory, similar to the abstract case. For instance, the profinite fundamental group of a graph of groups over a single edge is an amalgamated free product i.e.\ finite pushout (see \cite[Example 6.2.3(d)]{ribesgraph}), whilst the fundamental group over a ``cone graph" with trivial edge groups and cone-vertex group is a free product (see \cite[Example 6.2.3(b)]{ribesgraph}). However, as far as the author is aware, there has not been much work that studies infinite amalgamated free products of profinite groups.

Let us now see that in the context of Bass--Serre theory (at least), our notion of colimits internal to $\Pro$ gives the correct notion of infinite amalgamated free products of profinite groups. We remark that these are in general not infinite categorical pushouts in $\PGrp$: those do exist\footnote{Warning: In \cite{ribesamal} and \cite{venj}, the authors say that a pushout \emph{exists} if the canonical maps into it are injective, which is not what we mean here. In \cite{profinite} and \cite{ribesgraph}, the authors say that a pushout is \emph{proper} if the canonical maps are injective.} but may not be well-behaved (cf.\ \cite[Remark 2]{venj}). The reader should consult \cite{ribesgraph} for the terminology and facts that we are about to use.

Let $X\in\Pro$ and $\Gamma=C(X,*)$ be the profinite cone graph of $X$, which is connected (see \cite[Example 3.4.2]{ribesgraph}). That is, $\Gamma$ has vertices $V=X\sqcup\{*\}$ and edges $E=\overline{X}\cong X$, with incidence relations saying that an edge $\overline{x}\in E$ has source $*\in V$ and target $x\in V$. Suppose that $\G=(\G,\varpi,\Gamma)$ is a profinite graph of groups over $\Gamma$, where all edge groups, as well as the vertex group above $*$, are (isomorphic to) a fixed profinite group $H$, such that the restriction of $\partial_0$ to each edge group is an isomorphism of $H$. We shall identify all these copies of $H$. In particular, we can think of $H$ as a common subgroup of every vertex group $\G(x)$, $x\in V\backslash\{*\}$. Let us write $\theta_x\colon H\inj \G(x)$ for the restriction of $\partial_1$ to the edge group $\G(\overline{x})=H$. For ordinary groups (and the ordinary version of the graph of groups we just defined), the fundamental group is the amalgamated free product $*_H\, \G(x)$ (cf.\ \cite[Section I.4.4]{trees}). Let us think about what the profinite fundamental group of $\G$ looks like.

Note that $\Gamma$, being an inverse limit of finite trees, is a profinite tree and is simply connected (\cite[Proposition 3.3.3]{ribesgraph}). Moreover, the set of edges $E$ is a closed subspace of $\Gamma$. Therefore, as pointed out at the end of \cite[page 178]{ribesgraph}, the fundamental group $\Pi_1(\G)$ of $\G$ is rather easy to define: it is a profinite group $\Pi=\Pi_1(\G)\in\PGrp$ with the universal property of having a map $\beta\colon\G_V=\varpi^{-1}(V)\to\Pi$ in $\Pro$, such that the restriction of $\beta$ to each fibre is a group homomorphism $\beta_v\colon\G(v)\to\Pi$ (where $v\in V$), and such that $\beta_x\circ\theta_x=\beta_*$ for each $x\in X$.

We now construct an internal category $A=A_\Gamma$ associated to $\Gamma$. It has $A_0=X\sqcup\{*\}=V$ and $$A_1=\overline{X}\sqcup X\sqcup\{*\}=E\sqcup V,$$ with the obvious composition rule etc.: we think of $A_\Gamma$ as $\Gamma$ with an identity map added to each vertex. We also define an internal functor $P\in\Grp(\IFun(A_\Gamma,\Pro))$ as follows. It has $P_0=\G_V=\bigsqcup_{v\in V}\G(v)$, $p_0=\varpi|_{\G_V}$ and $$p_1=\partial_1\colon A_1\times_{A_0}\G_V=\G\to \G_V.$$ Note that $P(\overline{x})=\theta_x$ (where $\overline{x}\in E$) and that $P$ is in fact a group object in the internal functor category.

Consider the colimit of $P$, which exists by Theorem \ref{thm1}. By definition, it is a profinite group $C$ with the universal property of having a map $P\to\Delta_{A_\Gamma}(C)$ in $\Grp(\IFun(A_\Gamma,\Pro))$. It is easy to see that a map $P\to\Delta_{A_\Gamma}(C)$ is the same as a map $\alpha\colon P_0=\G_V\to C$ in $\Pro$ such that the restriction of $\alpha$ to each fibre is a group homomorphism $\alpha_v\colon\G(v)\to C$ (where $v\in V$) and such that all the $\alpha_v$ collectively define a natural transformation from $P$ to $\Delta_{A_\Gamma}(C)$, viewed as functors $A\to\Grp$. Moreover, to say that the $\alpha_v$ satisfy the naturality condition is the same as to say that for each $\overline{x}\in E=\overline{X}$, we have $\alpha_x\circ\theta_x=\alpha_*$. But we have just written down precisely the universal property of the fundamental group $\Pi_1(\G)$!

Thus, the fundamental group $\Pi_1(\G)$ of the specific graph of groups $\G$ over $\Gamma$ we defined is an infinite \emph{internal pushout} i.e.\ a colimit of shape $A_\Gamma$ internal to $\Pro$. We will denote this group by $\coprod_H\G(x)$ (not to be confused with the free profinite product) and should think of it as the \emph{free (profinite) product of the $\G(x)$ amalgamating $H$ along $\theta_x$}, where the topology on $\G$ is implicit. Since we have established $\coprod_H\G(x)$ as a fundamental group, we could use the results of \cite{ribesgraph} to study this group if we wish. As an example, the Mayer--Vietoris sequence (\cite[Theorem 9.4.1]{ribesgraph}) easily implies the following.

\begin{prop}
Let $G=\coprod_H\G(x)$ be as above and assume that it is proper, i.e.\ the canonical maps $\G(x)\to G$ are injective. Let $R$ be a profinite quotient ring of $\widehat{\Z}$ and $B$ be a profinite $\RG$-module which is projective when restricted to $H$. Moreover, write $H_i(G,B)$ for the $i$th homology group of $G$ with coefficients in $B$. Then, for every $i\geq2$, there is an isomorphism of profinite $R$-modules $$H_i(G,B)=\widehat{\bigoplus}_{x\in X}H_i(\G(x),B).\footnote{Recall that $\widehat{\oplus}$ refers to the coproduct of profinite $R$-modules internal to $\Pro$ (i.e.\ the ``profinite direct sum" in \cite{jc3}), analogous to what we have already done for groups. The symbol $\oplus$ (resp.\ $\boxplus$) is used in \cite{ribesgraph} (resp.\ \cite{gareth}) instead.}$$
\end{prop}

\begin{rmk}
The construction $\Gamma\mapsto A_\Gamma$ above might seem familiar. Recall that the forgetful functor $\Cat\to\Graph$, where $\Graph$ is the category of directed multigraphs, has a left adjoint which takes a graph $\Gamma$ and freely constructs a category $A_\Gamma$ from it (see \cite[Section II.7]{maccat}). Let $\PGraph$ be the category of profinite graphs, where by a profinite graph\footnote{This is less general than the profinite graphs appearing in \cite{ribesgraph}. In fact, the profinite graphs defined here are exactly those profinite graphs in \cite{ribesgraph} whose set of edges is closed.} here we mean a pair of profinite spaces $(V,E)$ together with continuous maps $d_0, d_1\colon E\to V$, and morphisms are required to send edges to edges. Then, there is clearly a forgetful functor $\Cat(\Pro)\to\PGraph$. Given a profinite graph $\Gamma=(V,E)\in\PGraph$, if every directed path in $\Gamma$ has length bounded by a fixed constant, then we can easily construct a category $A=A_\Gamma$ internal to $\Pro$ by following the argument of \cite[Section II.7]{maccat}. By a directed path, we just mean a sequence of composable arrows, possibly with repetition. Briefly, we have $A_0=V$ and $A_1=\bigsqcup^\infty_{n=0}E_n$, where $E_n$ is the space of directed paths in $\Gamma$ of length $n$. We need to bound the lengths of all paths in $\Gamma$ so that $A_1$ is in fact a finite disjoint union. It is then obvious how to make $A_\Gamma$ a category internal to $\Pro$ and that $A_\Gamma$ is the reflection of $\Gamma$ along the forgetful functor.

However, if we allow $\Gamma$ to have arbitrarily long paths, the theory quickly becomes tricky. This already happens when $\Gamma$ consists of just a single vertex and a single loop—then $A_\Gamma$ should be the one-object category corresponding to the ``free profinite monoid of rank 1"!
\end{rmk}

\begin{rmk}
Another important construction that appears in Bass--Serre theory is that of HNN-extensions. It is well known that these are not 1-colimits of groups, but 2-coinserters (a type of finite 2-colimits). We can develop a similar theory here, but since this is not the focus of the current writing, we will sketch the basic ideas below and leave it to the interested reader to fill the gaps.

The category $\Cat(\Pro)$ is naturally a 2-category if we take internal natural transformations as 2-morphisms. Observe that one-object internal categories are precisely profinite monoids, i.e.\ monoid objects internal to $\Pro$ (viewed as one-object categories), and internal functors between them are precisely continuous monoid homomorphisms. Thus, the 1-category $\PGrp$ can be viewed as a full 1-subcategory of $\Cat(\Pro)$ and hence $\PGrp$ inherits a 2-categorical structure. It is easy to see that given two profinite groups $G,H$ and two maps $\alpha,\beta\colon G\to H$, an internal natural transformation $t\colon \alpha\to \beta$ is the same as an element $t\in H$ such that $\alpha(g)=t^{-1}\beta(g)t$ for all $g\in G$. Thus, given $\alpha,\beta$ as above and assuming they are injective, the profinite HNN-extension of $\alpha$ and $\beta$ is exactly a profinite group $K$ with the universal property of having a 1-morphism $\gamma\colon H\to K$ and a 2-morphism $\gamma\circ\alpha\to\gamma\circ\beta$ i.e.\ it is a 2-coinserter.
\end{rmk}

\section{Left Adjoints and Coproducts}
\label{sec3}

Now we shall study how left adjoints and other functors interact with coproducts internal to $\Pro$. This will eventually lead to our main theorem (Theorem \ref{mainthm}) and main examples (Example \ref{maineg}). Let us once again forget about internal categories for a moment and think about $\Set$ and $\Grp$. The forgetful functor $\Grp\to\Set$ has a left adjoint $F$, namely the free groups functor. The fact that left adjoints preserve colimits in this situation can be expressed as follows: for any $A\in\Cat$, we have a commutative diagram, up to natural isomorphism, of (large) categories
\[\begin{tikzcd}
{\Fun(A,\Set)} \arrow[r, "\coprod_A"] \arrow[d, "F\circ(-)"'] & \Set \arrow[d, "F"] \\
{\Fun(A,\Grp)} \arrow[r, "\coprod_A"']                        & \Grp               
\end{tikzcd}\]

What would the corresponding diagram in the context of internal categories be? The top row should be $\IFun(A,\Pro)\overset{\coprod_A}{\longrightarrow}\Pro$, where $A\in\Cat(\Pro)$ and $\coprod_A$ is the left adjoint of $\Delta_A$, and similarly for the bottom row but for group objects. The left vertical map, however, presents difficulties, as it seems like we would need to view the free profinite groups functor $F\colon\Pro\to\PGrp$ as an internal functor. We can circumvent this problem by considering only coproducts rather than all colimits internal to $\Pro$.

Let $A=X^\delta$ be the discrete internal category associated to $X\in\Pro$. Then, we would like to have a commutative diagram, up to natural isomorphism, of ordinary categories
\[\begin{tikzcd}
\Pro_{/X} \arrow[r] \arrow[d, "F_X"']   & \Pro \arrow[d, "F"] \\
\Grp(\Pro_{/X}) \arrow[r, "\coprod_X"'] & \PGrp              
\end{tikzcd}
\tag{1}
\label{cd1}\]
where $F$ is any functor (not necessarily a left adjoint for now) and the top arrow simply sends an object $(P,X)=(P\to X)$ to $P$. Next, we need to construct the functor $F_X$, which should be induced by $F$. A first idea is that for an object $p\colon P\to X$ in $\Pro_{/X}$, the object $F_X(P)$ should have fibres $F_X(P)(x)=F(P(x))$, $x\in X$, where $P(x)=p^{-1}(\{x\})$. The issue is that we still need to give \begin{equation}\label{eqn1} F_X(P)=\bigsqcup_{x\in X} F(P(x))\tag{$\star$}\end{equation} an appropriate profinite topology. We will do so by using the fact that all the categories above are subcategories of certain pro-completions (refer to \cite{elephant} for the necessary background). What follows should be viewed as a more general version of \cite[Proposition A.13]{garethrel}.

\subsection{Preservation of Coproducts Internal to $\Pro$}\label{3.1}

Let us first set up appropriate notions for the three main types of categories that we will be working with. The definitions below might seem sporadic and one might ask if there is a uniform way to describe them. We will address this question in Remark \ref{general}.

\begin{itemize}
\item Let $\BPro=\Fun(\mathbf{2},\Pro)$ be the category of bundles of profinite spaces, where $\mathbf{2}$ is the interval category. That is, $\BPro$ has objects continuous maps $P\to X$ of profinite spaces, and morphisms commutative squares. Let $\BProf$ be the full subcategory of bundles $P\to X$ where both $P$ and $X$ are finite. It is a classical fact (see \cite[Proposition 8.8.5]{gro1}) that $\BPro=\Pro(\BProf).$\footnote{Here, $\Pro(\mathbf{C})$ denotes the pro-completion of $\mathbf{C}$. Recall that we also write $\Pro$ for the category $\Pro(\mathbf{Set}_\mathrm{fin})$ of profinite spaces, which hopefully doesn't cause any confusion.}

\item Let $\BPGrp$ be the category of bundles of profinite groups (see Remark \ref{rmk1}). That is, $\BPGrp$ has objects pairs $(P,X)$ with $X\in\Pro$ and $P\in\Grp(\Pro_{/X})$. Let $\BPGrpf$ be the full subcategory of bundles where both spaces are finite. It follows from the group analogue of \cite[Corollary 2.13]{jc3} that $\BPGrp=\Pro(\BPGrpf)$.

\item Fix a profinite ring $R$. Let $\BPModR$\footnote{This category is (equivalent to the one) denoted by $\coshR$ in \cite{jc3}.} be the category of bundles of profinite $R$-modules (see \cite[Section 9.1]{ribesgraph} or Remark \ref{rmk2}). That is, $\BPModR$ has objects pairs $(P,X)$ with $X\in\Pro$ and $P$ an $(R\times X)$-module object in $\Pro_{/X}$. Let $\BPModRf$ be the full subcategory of bundles where both spaces are finite. It follows from \cite[Corollary 2.13]{jc3} that $\BPModR=\Pro(\BPModRf)$.
\end{itemize}

Let $\mathbf{C}$ denote one of $\mathbf{Set}_\mathrm{fin}$ (the category of finite sets), $\mathbf{Grp}_\mathrm{fin}$ (the category of finite groups), and $\mathbf{Mod}(R)_\mathrm{fin}$ (the category of finite topological $R$-modules for a fixed profinite ring $R$), and let $\mathbf{BC}$ denote $\BProf$, $\BPGrpf$, and $\BPModRf$ respectively. Observe that an object $(P,X)\in\mathbf{BC}$ is equivalently a finite family $\{P(x)\}_{x\in X}$ of objects in $\mathbf{C}$. Also, $\Pro(\mathbf{BC})$ will be equivalent to $\BPro$, $\BPGrp$, and $\BPModR$ respectively, a fact which can be conveniently expressed by writing $\Pro(\mathbf{BC})=\mathbf{B}\Pro(\mathbf{C})$.

For a fixed profinite space $X$ and $\mathbf{C}=\mathbf{Set}_\mathrm{fin}$ (resp.\ $\mathbf{Grp}_\mathrm{fin}$, resp.\ $\mathbf{Mod}(R)_\mathrm{fin}$), let $\Pro(\mathbf{C})_X$ denote the category $\Pro_{/X}$ (resp.\ $\Grp(\Pro_{/X})$, resp.\ the category of $(R\times X)$-module objects in $\Pro_{/X}$)\footnote{In other words, $\Pro(\mathbf{C})_X$ is the fibre of the Grothendieck construction $\mathbf{B}\Pro(\mathbf{C})$ at $X\in\Pro$.}.

To state our main theorem, we still need the definition of \emph{relative adjunctions} from \cite{ulmer}. We will look at some examples of relative adjunctions in Section \ref{3.2}.

\begin{defn}
Let $\mathbf{E}$, $\mathbf{F}$, $\mathbf{G}$ be ordinary categories and $L\colon \mathbf{E}\to \mathbf{G}$, $R\colon \mathbf{F}\to \mathbf{E}$, $J\colon \mathbf{F}\to \mathbf{G}$ functors. We say that $L$ is \emph{$J$-left coadjoint to R} if for every $e\in\mathbf{E}$ and $f\in\mathbf{F}$, there is a natural isomorphism $$\Hom_{\mathbf{G}}(Le,Jf)=\Hom_{\mathbf{E}}(e,Rf).$$
\end{defn}

The following is then our main theorem. In Section \ref{3.3}, we will see many examples that fit into the following theorem.

\begin{thm}\label{mainthm}
Let $\mathbf{C}$ and $\mathbf{D}$ each be any of $\Set_\mathrm{fin}$, $\Grp_\mathrm{fin}$, or $\ModR_\mathrm{fin}$, where $R$ is a profinite ring\footnote{We allow $\mathbf{C}=\ModR_\mathrm{fin}$ and $\mathbf{D}=\mathbf{Mod}(S)_\mathrm{fin}$ for different profinite rings $R$ and $S$.}. Let $J\colon\mathbf{D}\inj\Pro(\mathbf{D})$ and $J^*\colon\mathbf{BD}\inj\Pro(\mathbf{BD})$ be the canonical embeddings.

Let $L\colon\mathbf{C}\to\Pro(\mathbf{D})$ be any functor and let $L'\colon\Pro(\mathbf{C})\to\Pro(\mathbf{D})$ be its extension. Then:

\begin{enumerate}[label=(\roman*)]
\item\label{mainthm1} The functor $L$ induces a functor $L^*\colon\mathbf{BC}\to\Pro(\mathbf{BD})=\mathbf{B}\Pro(\mathbf{D})$ by sending a finite bundle $(P,X)=\bigsqcup_{x\in X}P(x)\in\mathbf{BC}$ to the disjoint union $\bigsqcup_{x\in X}LP(x)\in\mathbf{B}\Pro(\mathbf{D})$, which then extends to a functor ${L^*}'\colon\mathbf{B}\Pro(\mathbf{C})\to\mathbf{B}\Pro(\mathbf{D})$. Suppose further that $L'$ commutes with finite coproducts. Then, we obtain a commutative diagram (up to natural isomorphism)
\[
\begin{tikzcd}
\mathbf{B}\Pro(\mathbf{C}) \arrow[d, "{L^*}'"'] \arrow[r, "\coprod"] & \Pro(\mathbf{C}) \arrow[d, "L'"] \\
\mathbf{B}\Pro(\mathbf{D}) \arrow[r, "\coprod"']                     & \Pro(\mathbf{D})                
\end{tikzcd}\]
Moreover, the restriction of the above diagram gives a well-defined commutative diagram
\[\begin{tikzcd}
\Pro(\mathbf{C})_{X} \arrow[d, "{L^*}'"'] \arrow[r, "\coprod_X"] & \Pro(\mathbf{C}) \arrow[d, "L'"] \\
\Pro(\mathbf{D})_{X} \arrow[r, "\coprod_X"']                     & \Pro(\mathbf{D})            
\end{tikzcd}\]

\item\label{mainthm2} Suppose $L'$ arises from Lemma \ref{lem1}. That is, suppose $L$ is $J$-left coadjoint to a functor $R\colon\mathbf{D}\to\mathbf{C}$, so $L'$ is left adjoint to $\Pro(R)$ and automatically commutes with finite coproducts. Let $R^*\colon\mathbf{BD}\to\mathbf{BC}$ be constructed from $R$ in the same way as $L^*$ from $L$. Then, $L^*$ is $J^*$-left coadjoint to $R^*$, ${L^*}'$ is left adjoint to $\Pro(R^*)$, and all eight functors appearing in the two diagrams above are left adjoints.
\end{enumerate}
\begin{proof}
Part \ref{mainthm1} will be proved in Propositions \ref{commprop} and \ref{bigthm}, and part \ref{mainthm2} will be proved in Lemma \ref{adjlemma} and Proposition \ref{leftadjgood}.
\end{proof}
\end{thm}

\begin{rmk}\label{general}
In fact, we can allow more general algebraic theories with coproducts internal to $\Pro$ and with the property that ``bundles commute with pro-completions". To make this precise, let $\mathbb{T}$ be a (finitary) Lawvere theory and write $\mathbf{Alg}_{\mathbb{T}}(\mathbf{S})$ for the category of $\mathbb{T}$-algebras in $\mathbf{S}$, where $\mathbf{S}$ is a category with finite products. For background on Lawvere theories, the reader may refer to \cite{lawvere}.

Let us write $\mathbf{B}\mathbf{Alg}_{\mathbb{T}}(\Pro)$ for the category whose objects are pairs $(P,X)$ with $X\in\Pro$ and $P\in\mathbf{Alg}_{\mathbb{T}}(\Pro_{/X})$. A morphism from $(P,X)$ to $(Q,Y)$ consists of continuous maps $f\colon P\to Q$ and $X\to Y$ making the obvious square commute, such that the restriction of $f$ to each fibre above $X$ is a $\mathbb{T}$-map. Observe that $\mathbf{B}\mathbf{Alg}_{\mathbb{T}}(\Pro)$ has all limits which are computed on fibres. Let $\mathbf{B}\mathbf{Alg}_{\mathbb{T}}(\Set_{\mathrm{fin}})$ be the full subcategory of $\mathbf{B}\mathbf{Alg}_{\mathbb{T}}(\Pro)$ where both spaces $P$ and $X$ are finite.

One can check that we only need the following two conditions on $\mathbb{T}$ for Theorem \ref{mainthm} to work for $\mathbf{C}=\mathbf{Alg}_{\mathbb{T}}(\Set_{\mathrm{fin}})$.
\begin{enumerate}[label=(\alph*)]
\item For every $X\in\Pro$, the diagonal functor $\Delta_{X^\delta}\colon\mathbf{Alg}_{\mathbb{T}}(\Pro)\to\mathbf{Alg}_{\mathbb{T}}(\Pro_{/X})$ has a left adjoint, i.e.\ the theory $\mathbb{T}$ has all coproducts internal to $\Pro$. In particular, $\mathbf{Alg}_{\mathbb{T}}(\Pro)$ has ordinary finite coproducts.
\item The canonical map $\Pro(\mathbf{B}\mathbf{Alg}_{\mathbb{T}}(\Set_{\mathrm{fin}}))\to\mathbf{B}\mathbf{Alg}_{\mathbb{T}}(\Pro)$ is an equivalence. In particular, the canonical map $\Pro(\mathbf{Alg}_{\mathbb{T}}(\Set_{\mathrm{fin}}))\to\mathbf{Alg}_{\mathbb{T}}(\Pro)$ is an equivalence, i.e.\ pro-(finite $\mathbb{T}$-algebras) are exactly $\mathbb{T}$-algebras internal to $\Pro$.
\end{enumerate}

We then recover $\mathbf{C}=\mathbf{Set}_{\mathrm{fin}}$ (resp.\ $\mathbf{Grp}_{\mathrm{fin}}$) by taking $\mathbb{T}$ to be the theory of sets (resp.\ groups). The problem with this is that the standard single-sorted Lawvere theory of $R$-modules does not carry any topology on $R$, so we do not obtain the theory of profinite $R$-modules from this (for a fixed profinite ring $R$). One solution is to use enriched Lawvere theories, but we immediately run into the problem that $\Pro$ is not Cartesian closed. We have some ideas on how to resolve this, but that would distract us from our current goal, so we have decided to list the categories we will work with, rather than to describe them in a uniform way.
\end{rmk}

Before we proceed, let us note some properties of $\mathbf{C}=\Set_\mathrm{fin}$, $\Grp_{\mathrm{fin}}$, or $\ModR_\mathrm{fin}$ and the associated bundle category. Firstly, inverse limits in $\mathbf{B}\Pro(\mathbf{C})$ commute with taking fibres (cf.\ \cite[Proposition 2.7]{jc3}). More precisely, if $$(P,X)=\varprojlim(P_i,X_i)=(\varprojlim P_i,\varprojlim X_i)$$ and $x=(x_i)\in\varprojlim X_i$, then we have $P(x)=\varprojlim P_i(x_i)$. We will sometimes write $P$ instead of $(P,X)$ when $X$ is understood. Furthermore, there is an internal coproduct functor $$\coprod\colon\mathbf{B}\Pro(\mathbf{C})\to\Pro(\mathbf{C})$$ which is left adjoint to the functor $Z\mapsto (Z,*)$. More accurately, if $\mathbf{C}=\mathbf{Set}_\mathrm{fin}$ (resp.\ $\mathbf{Grp}_\mathrm{fin}$, resp.\ $\mathbf{Mod}(R)_\mathrm{fin}$), then $\coprod$ sends a bundle $(P,X)$ to $P$ (resp.\ the free product $\coprod_XP(x)$, resp.\ the profinite direct sum $\widehat{\oplus}_XP(x)$), but for ease of notation, we will use the coproduct symbol $\coprod$ or $\coprod_X$ to denote all of them.

Suppose each of $\mathbf{C}$ and $\mathbf{D}$ is a category listed above, and $L\colon\mathbf{C}\to\Pro(\mathbf{D})$ is any functor. By the universal property of pro-completions, $L$ extends uniquely to a functor $L'\colon\Pro(\mathbf{C})\to\Pro(\mathbf{D})$ which commutes with inverse limits.

There is an induced functor $L^*\colon\mathbf{BC}\to\mathbf{B}\Pro(\mathbf{D})$ which sends a finite bundle $P\to X$ to the bundle $$\bigsqcup_{x\in X}L(P(x)) \to X,$$ where the finite disjoint union is given the disjoint union topology, and sends morphisms in the obvious way. Note that this is well defined, as for each $x\in X$, the fibre $P(x)$ is in $\mathbf{C}$. By the universal property of pro-completions again, $L^*$ uniquely extends to a functor $${L^*}'\colon\Pro(\mathbf{BC})=\mathbf{B}\Pro(\mathbf{C})\to\mathbf{B}\Pro(\mathbf{D}).$$ By construction, for $(P,X)=\varprojlim(P_i,X_i)\in\Pro(\mathbf{BC})$ with $(P_i,X_i)\in\mathbf{BC}$, we have that ${L^*}'(P,X)=\varprojlim L^*(P_i,X_i),$ which is a bundle over the same space $X$. If moreover $x=(x_i)\in X$, then $${L^*}'(P)(x)=(\varprojlim L^*P_i)(x)=\varprojlim((L^*P_i)(x_i))=\varprojlim L(P_i(x_i))=L'(\varprojlim P_i(x_i))=L'(P(x)).$$ This is precisely our initial motivation, except ${L^*}'(P)=\varprojlim L^*(P_i)$ is also canonically equipped with the inverse limit topology: see (\ref{eqn1}).

We then obtain the following analogue of \cite[Proposition A.13]{garethrel}.

\begin{prop}\label{commprop}
Suppose the functor $L'\colon\Pro(\mathbf{C})\to\Pro(\mathbf{D})$ commutes with finite coproducts. Then, the following diagram commutes (up to natural isomorphism):
\[
\begin{tikzcd}
\mathbf{B}\Pro(\mathbf{C}) \arrow[r, "\coprod"] \arrow[d, "{L^*}'"'] & \Pro(\mathbf{C}) \arrow[d, "L'"] \\
\mathbf{B}\Pro(\mathbf{D}) \arrow[r, "\coprod"']      & \Pro(\mathbf{D})              
\end{tikzcd}
\tag{2}\label{cd2}
\]
\begin{proof}
Note that all the functors in the above diagram commute with inverse limits: $L'$ and ${L^*}'$ do by definition, and for the other two see \cite[Proposition 5.1.7]{ribesgraph} or \cite[Proposition 2.7]{jc3}. By the universal property of pro-completions, it suffices to check that the diagram commutes when restricted to $\mathbf{BC}$. In this case, the top-right composition sends a finite bundle $(P,X)$ to $L'(\coprod_{x\in X}P(x))$ and the bottom-left composition sends it to $$\coprod_{x\in X}L(P(x))=\coprod_{x\in X}L'(P(x)).$$ These are the same precisely because $L'$ commutes with finite coproducts. We remind the reader that finite coproducts internal to $\Pro$ are just ordinary categorical coproducts (see Example \ref{egdisc}).
\end{proof}
\end{prop}

\begin{rmk}\label{what}
Suppose $U\colon\Pro(\mathbf{C})\to\Pro(\mathbf{D})$ is any right adjoint. Then, $U$ preserves inverse limits, so must be the unique extension of its restriction to $\mathbf{C}$. Thus, the diagram
\[\begin{tikzcd}
\mathbf{B}\Pro(\mathbf{C}) \arrow[r, "\coprod"]   \arrow[d, "{U^*}'"']                       & \Pro(\mathbf{C}) \arrow[d, "U"]                      \\
\mathbf{B}\Pro(\mathbf{D}) \arrow[r, "\coprod"'] & \Pro(\mathbf{D}) 
\end{tikzcd}\]
commutes if $U$ commutes with finite coproducts. If $\mathbf{C}=\ModR_\mathrm{fin}$ and $\mathbf{D}=\mathbf{Mod}(S)_\mathrm{fin}$, so that finite products and coproducts agree in $\PModR$ and in $\mathbf{PMod}(S)$ (where $\PModR$ is the category of profinite $R$-modules), then we get the somewhat unintuitive statement that right adjoints \emph{always} preserve (internal) coproducts of profinite modules. This should sound a lot more reasonable once we realise that the internal coproduct $\widehat{\oplus}$ of profinite modules, despite being defined as a left adjoint, actually generalises categorical \emph{products} of modules (see \cite[Example 5.6.4(c)]{ribesgraph} or Example \ref{maineg}\ref{maineg4} below).
\end{rmk}

Now, let us see how we can obtain a commutative diagram similar to (\ref{cd1}). Note that $\Pro(\mathbf{C})_X$ is a non-full subcategory of $\mathbf{B}\Pro(\mathbf{C})$.

\begin{prop}\label{bigthm}
For each $X\in\Pro$, the restriction of diagram (\ref{cd2}) gives a well-defined diagram
\[
\begin{tikzcd}
\Pro(\mathbf{C})_X \arrow[r, "\coprod_X"] \arrow[d, "{L^*}'"'] & \Pro(\mathbf{C}) \arrow[d, "L'"] \\
\Pro(\mathbf{D})_X \arrow[r, "\coprod_X"']      & \Pro(\mathbf{D})              
\end{tikzcd}
\tag{3}\label{cd3}
\]
which commutes if (\ref{cd2}) does.
\begin{proof}
As mentioned in the previous paragraph, the only problem is that $\Pro(\mathbf{D})_X$ is a non-full subcategory of $\mathbf{B}\Pro(\mathbf{D})$, so we need to ensure that ${L^*}'$ sends objects and morphisms in $\Pro(\mathbf{C})_X$ to those in $\Pro(\mathbf{D})_X$.

Given a morphism $f\colon P\to Q$ in $\Pro(\mathbf{C})_X$, decompose it into an inverse limit in $\mathbf{B}\Pro(\mathbf{C})=\Pro(\mathbf{BC})$
\[\begin{tikzcd}
P_i \arrow[rr, "f_i"] \arrow[rd] &     & Q_i \arrow[ld] \\
                                 & X_i &               
\end{tikzcd}\]
where all spaces are finite. By construction, the morphism ${L^*}'(f)$ is given by the inverse limit of
\[\begin{tikzcd}
L^*(P_i) \arrow[rr, "L^*(f_i)"] \arrow[rd] &     & L^*(Q_i) \arrow[ld] \\
                                           & X_i &                    
\end{tikzcd}\]
where $L^*(P_i)\to X_i$ is the obvious map $\bigsqcup_{x\in X_i}L(P_i(x)) \to X_i$ and similarly for $Q_i$. It is thus clear that the inverse limit of the above diagrams in fact lives in $\Pro(\mathbf{D})_X$.
\end{proof}
\end{prop}

\subsection{Relative Adjunctions}\label{3.2} In Proposition \ref{commprop}, we required $L'$, the extension of $$L\colon\mathbf{C}\to\Pro(\mathbf{D}),$$ to commute with finite coproducts, not that it be a left adjoint. However, many such functors $L'$ which occur naturally in the profinite setting also happen to be left adjoints (see Examples \ref{exampleJco}\ref{exampleJco1} and \ref{maineg}). In this case, we can say more about diagrams (\ref{cd2}) and (\ref{cd3}). To explain why left adjoints are common in this setting, we need the concept of \emph{relative adjoints}.

Let $\mathbf{E}$, $\mathbf{F}$, $\mathbf{G}$ be ordinary categories and $L\colon \mathbf{E}\to \mathbf{G}$, $R\colon \mathbf{F}\to \mathbf{E}$, $J\colon \mathbf{F}\to \mathbf{G}$ functors. Recall from (the dual of) \cite[Definition 2.2]{ulmer} that $L$ is \emph{$J$-left coadjoint to R} if for every $e\in\mathbf{E}$ and $f\in\mathbf{F}$, there is a natural isomorphism $$\Hom_{\mathbf{G}}(Le,Jf)=\Hom_{\mathbf{E}}(e,Rf).$$ If $J$ is the identity functor, then we recover the usual definition of a left adjoint. The following simple observation is why relative adjoints appear often in the study of profinite groups. In the following lemma only, $\mathbf{C}$ and $\mathbf{D}$ are not necessarily $\mathbf{Set}_\mathrm{fin}$, $\mathbf{Grp}_\mathrm{fin}$, or $\mathbf{Mod}(R)_\mathrm{fin}$.

\begin{lemma}(cf.\ \cite{proadj})\label{lem1}
Let $\mathbf{C}$ and $\mathbf{D}$ be small categories with finite limits, $J\colon\mathbf{D}\inj\Pro(\mathbf{D})$ be the canonical embedding, and $L\colon\mathbf{C}\to\Pro(\mathbf{D})$, $R\colon\mathbf{D}\to\mathbf{C}$ be functors such that $L$ is $J$-left coadjoint to $R$. Let $$L'\colon\Pro(\mathbf{C})\to\Pro(\mathbf{D})\,\,\,\, \text{and}\,\,\,\, \Pro(R)\colon\Pro(\mathbf{D})\to\Pro(\mathbf{C})$$ be the canonical extensions. Then, $L'$ is left adjoint to $\Pro(R)$.
\begin{proof}
Let $c=\varprojlim c_i\in\Pro(\mathbf{C})$ and $d=\varprojlim d_j\in\Pro(\mathbf{D})$. Then, we have natural isomorphisms
\begin{eqnarray*}
\Hom_{\Pro(\mathbf{D})}(L'c,d)&=&\varprojlim_j\varinjlim_i\Hom_{\Pro(\mathbf{D})}(Lc_i,Jd_j)  \\
&=&\varprojlim_j\varinjlim_i\Hom_{\mathbf{C}}(c_i,Rd_j)  \\
&=&\Hom_{\Pro(\mathbf{C})}(c,\varprojlim Rd_j)  \\
&=&\Hom_{\Pro(\mathbf{C})}(c,\Pro(R)d),
\end{eqnarray*}
which completes the proof.
\end{proof}
\end{lemma}

\begin{exmp}\label{exampleJco}
\begin{enumerate}[label=(\roman*)]
\item\label{exampleJco1} Take $\mathbf{C}=\Set_{\mathrm{fin}}$ and $\mathbf{D}=\Grp_{\mathrm{fin}}$. Let $L\colon\Set_{\mathrm{fin}}\to\PGrp$ be the free profinite groups functor, i.e.\ $L(X)$ is the free profinite group on the finite set $X$. This can be constructed as the profinite completion of the abstract free group on $X$. Let $R\colon\Grp_{\mathrm{fin}}\to\Set_{\mathrm{fin}}$ be the forgetful functor. It is immediate that $L$ is $J$-left coadjoint to $R$ (where $J\colon\Grp_{\mathrm{fin}}\inj\PGrp$), so by the above lemma, we get the free profinite groups functor $L'\colon\Pro\to\PGrp$ which is left adjoint to the forgetful functor $\PGrp\to\Pro$. We used here the fact that the forgetful functor $\PGrp\to\Pro$ commutes with inverse limits, so that it is indeed $\Pro(R)$. By construction, we have $$L'(\varprojlim X_i)=\varprojlim L(X_i)$$ for $X_i$ finite. The point is that it's easy to construct $L$ and check that it is $J$-left coadjoint to $R$ because only finite bases are involved, and then we automatically get an actual adjunction between the pro-completions.

This works in various other situations in the profinite setting. To name two, notice that the free modules functor $$R\llbracket-\rrbracket\colon\Pro\to\PModR$$ and induction $$\Ind^G_H(-)\colon\mathbf{PMod}(\kH)\to\mathbf{PMod}(\kG)$$ both arise this way.
\item\label{exampleJco2} Let $\mathbf{C}$ be one of $\mathbf{Set}_\mathrm{fin}$, $\mathbf{Grp}_\mathrm{fin}$, and $\mathbf{Mod}(R)_\mathrm{fin}$. In fact, the internal coproduct functor $\coprod\colon\mathbf{B}\Pro(\mathbf{C})\to\Pro(\mathbf{C})$ also arises from Lemma \ref{lem1}: it is the extension of $$\coprod\colon\mathbf{BC}\to\Pro(\mathbf{C}),$$ which is $J$-left coadjoint to the functor $$\mathbf{C}\to\mathbf{BC}\colon Z\mapsto(Z,*),$$ where $J\colon\mathbf{C}\inj\Pro(\mathbf{C})$. Once again, because $\mathbf{BC}$ consists of finite bundles $(P,X)$, the internal coproduct functor sends $(P,X)\in\mathbf{BC}$ to the ordinary coproduct $\coprod_{x\in X}P(x)$, so left coadjointness is immediate.
\end{enumerate}
\end{exmp}

In the setting of bundles in Section \ref{3.1}, suppose $L'\colon \Pro(\mathbf{C})\to\Pro(\mathbf{D})$ is actually a left adjoint which arises from Lemma \ref{lem1}. That is, suppose $L\colon\mathbf{C}\to\Pro(\mathbf{D})$ is $J$-left coadjoint to a functor $R\colon\mathbf{D}\to\mathbf{C},$ where $J\colon\mathbf{D}\inj\Pro(\mathbf{D})$ is the canonical embedding. Let $$J^*\colon \mathbf{BD}\inj\Pro(\mathbf{BD})=\mathbf{B}\Pro(\mathbf{D})$$ be the canonical embedding and $R^*\colon\mathbf{BD}\to\mathbf{BC}$ be constructed from $R$ in the same way as $L^*$ from $L$. That is, $R^*$ sends a finite bundle $Q\to Y$ in $\mathbf{BD}$ to the bundle $$\bigsqcup_{y\in Y}R(Q(y))\to Y.$$ Then, it's easy to see that $L^*\colon\mathbf{BC}\to\mathbf{B}\Pro(\mathbf{D})$ is $J^*$-left coadjoint to $R^*\colon\mathbf{BD}\to\mathbf{BC}$. By Lemma \ref{lem1}, we obtain the following.

\begin{lemma}\label{adjlemma}
Suppose $L'\colon\Pro(\mathbf{C})\to\Pro(\mathbf{D})$ is a left adjoint that arises from $L\colon\mathbf{C}\to\Pro(\mathbf{D})$ being $J$-left coadjoint to $R\colon\mathbf{D}\to\mathbf{C}$. Then, the extension ${L^*}'\colon\mathbf{B}\Pro(\mathbf{C})\to\mathbf{B}\Pro(\mathbf{D})$ is left adjoint to $\Pro(R^*)\colon\mathbf{B}\Pro(\mathbf{D})\to\mathbf{B}\Pro(\mathbf{C})$.
\end{lemma}

Moreover, in this case $L'$ automatically commutes with finite coproducts because it is a left adjoint. By Proposition \ref{commprop}, we obtain the following commutative diagram consisting of four left adjoints:
\[
\begin{tikzcd}
\mathbf{B}\Pro(\mathbf{C}) \arrow[r, "\coprod"] \arrow[d, "{L^*}'"'] & \Pro(\mathbf{C}) \arrow[d, "L'"] \\
\mathbf{B}\Pro(\mathbf{D}) \arrow[r, "\coprod"']      & \Pro(\mathbf{D})              
\end{tikzcd}
\tag{2}\label{cd2'}
\]
The four right adjoints (cf.\ Example \ref{exampleJco}\ref{exampleJco2}) are given by
\[
\begin{tikzcd}
\mathbf{B}\Pro(\mathbf{C})                        & \Pro(\mathbf{C}) \arrow[l, "{Z\mapsto(Z,*)}"']                      \\
\mathbf{B}\Pro(\mathbf{D}) \arrow[u, "\Pro(R^*)"] & \Pro(\mathbf{D}) \arrow[l, "{Z\mapsto(Z,*)}"] \arrow[u, "\Pro(R)"']
\end{tikzcd}
\]

We already knew that (\ref{cd2'}) is commutative, but this also follows from the easy observation that the diagram just above is the pro-completion of the commutative diagram
\[\begin{tikzcd}
\mathbf{BC}           & \mathbf{C} \arrow[l,"{Z\mapsto(Z,*)}"']           \\
\mathbf{BD}\arrow[u,"R^*"] & \mathbf{D} \arrow[l,"{Z\mapsto(Z,*)}"] \arrow[u,"R"']
\end{tikzcd}\]
Thus, in situations more general than ours, it is often still easy to check whether a diagram like (\ref{cd2'}) commutes as long as all four functors arise from Lemma \ref{lem1}. Let us make this a formal observation.

\begin{lemma}\label{lem2}
Let $\mathbf{C}_1, \mathbf{C}_2, \mathbf{C}_3, \mathbf{C}_4$ be small categories with finite limits and let $J_i\colon\mathbf{C}_i\inj\Pro(\mathbf{C}_i)$ be the canonical embeddings. Suppose that there are some functors $R_{ij}\colon\mathbf{C}_i\to\mathbf{C}_j$ making the following square commute
\[\begin{tikzcd}
\mathbf{C}_1                     & \mathbf{C}_2 \arrow[l, "R_{21}"']                     \\
\mathbf{C}_3 \arrow[u, "R_{31}"] & \mathbf{C}_4 \arrow[l, "R_{43}"] \arrow[u, "R_{42}"']
\end{tikzcd}\]
and that, for each $(i,j)$ appearing above, there is a functor $L_{ji}\colon\mathbf{C}_j\to\Pro(\mathbf{C}_i)$ which is $J_i$-left coadjoint to $R_{ij}$. Let $L_{ji}'\colon\Pro(\mathbf{C}_j)\to\Pro(\mathbf{C}_i)$ be the (unique) extension of $L_{ji}$. Then, the following squares commute, where each map in the left square is left adjoint to the corresponding map in the right square.
\[\begin{tikzcd}
\Pro(\mathbf{C}_1) \arrow[r, "L_{12}'"] \arrow[d, "L_{13}'"'] & \Pro(\mathbf{C}_2) \arrow[d, "L_{24}'"] &  & \Pro(\mathbf{C}_1)                           & \Pro(\mathbf{C}_2) \arrow[l, "\Pro(R_{21})"']                           \\
\Pro(\mathbf{C}_3) \arrow[r, "L_{34}'"']                      & \Pro(\mathbf{C}_4)                      &  & \Pro(\mathbf{C}_3) \arrow[u, "\Pro(R_{31})"] & \Pro(\mathbf{C}_4) \arrow[l, "\Pro(R_{43})"] \arrow[u, "\Pro(R_{42})"']
\end{tikzcd}\]
\end{lemma}

As in Proposition \ref{bigthm}, we can restrict diagram (\ref{cd2'}) to get a commutative diagram (\ref{cd3}) with all the functors being left adjoints, but we have to be careful as we are dealing with non-full subcategories.

\begin{prop}\label{leftadjgood}
Suppose $L'\colon\Pro(\mathbf{C})\to\Pro(\mathbf{D})$ is a left adjoint that arises from $L\colon\mathbf{C}\to\Pro(\mathbf{D})$ being $J$-left coadjoint to $R\colon\mathbf{D}\to\mathbf{C}$, so we have a commutative diagram
\[
\begin{tikzcd}
\mathbf{B}\Pro(\mathbf{C}) \arrow[r, "\coprod"] \arrow[d, "{L^*}'"'] & \Pro(\mathbf{C}) \arrow[d, "L'"] \\
\mathbf{B}\Pro(\mathbf{D}) \arrow[r, "\coprod"']      & \Pro(\mathbf{D})              
\end{tikzcd}
\tag{2}\label{cd2''}
\]
Then, for each $X\in\Pro$, we can restrict the above to give a well-defined commutative diagram
\[
\begin{tikzcd}
\Pro(\mathbf{C})_X \arrow[r, "\coprod_X"] \arrow[d, "{L^*}'"'] & \Pro(\mathbf{C}) \arrow[d, "L'"] \\
\Pro(\mathbf{D})_X \arrow[r, "\coprod_X"']      & \Pro(\mathbf{D})              
\end{tikzcd}
\tag{3}\label{cd3''}
\]
where all four functors are left adjoints. The four right adjoints are given by
\[
\begin{tikzcd}
\Pro(\mathbf{C})_X                        & \Pro(\mathbf{C}) \arrow[l, "Z\mapsto Z\times X"']                      \\
\Pro(\mathbf{D})_X \arrow[u, "\Pro(R^*)"] & \Pro(\mathbf{D}) \arrow[l, "Z\mapsto Z\times X"] \arrow[u, "\Pro(R)"']
\end{tikzcd}
\]
\begin{proof}
In view of Proposition \ref{bigthm}, we only need to check that the right adjoint $$\Pro(R^*)\colon\mathbf{B}\Pro(\mathbf{D})\to\mathbf{B}\Pro(\mathbf{C})$$ of ${L^*}'$ restricts correctly and is still right adjoint. The exact same argument of Proposition \ref{bigthm} shows that $\Pro(R^*)$ restricts to the subcategories correctly.

Write $V=\Pro(R^*)$ for simplicity. To show that the restrictions ${L^*}'\colon\Pro(\mathbf{C})_X\to\Pro(\mathbf{D})_X$ and $V\colon\Pro(\mathbf{D})_X\to\Pro(\mathbf{C})_X$ are still adjoints, we need to check that given $P\in\Pro(\mathbf{C})_X$, $Q\in\Pro(\mathbf{D})_X$ and a map $P\to V(Q)$ in $\Pro(\mathbf{C})_X$, the unique induced map ${L^*}'(P)\to Q$ in $\mathbf{B}\Pro(\mathbf{D})$ in fact lives in $\Pro(\mathbf{D})_X$ (and vice versa, which is the same argument).

Let $\epsilon\colon {L^*}'\circ V\to \mathrm{id}$ be the counit of the adjunction between $\mathbf{B}\Pro(\mathbf{C})$ and $\mathbf{B}\Pro(\mathbf{D})$, so given $Q\in\mathbf{B}\Pro(\mathbf{D})$, the component $\epsilon_Q$ is the image of the identity map under $$\Hom_{\mathbf{B}\Pro(\mathbf{C})}(V(Q),V(Q))\to\Hom_{\mathbf{B}\Pro(\mathbf{D})}({L^*}'V(Q),Q).$$ Observe that if $Q$ is in $\Pro(\mathbf{D})_X$, then $\epsilon_Q$ in fact lives in $\Hom_{\Pro(\mathbf{D})_X}({L^*}'V(Q),Q)$ (cf.\ the isomorphisms in Lemma \ref{lem1}), which completes the proof.
\end{proof}
\end{prop}

For the rest of the paper, we will write $\PGrp_X$ (resp.\ $\PModR_X$) for $\Pro(\mathbf{C})_X$ when $\mathbf{C}=\Grp_{\mathrm{fin}}$ (resp.\ $\ModR_{\mathrm{fin}}$). In particular, we will write $\PAb_X=\mathbf{PMod}(\widehat{\Z})_X$ for the category of bundles of profinite abelian groups over $X\in\Pro$.

\begin{rmk}\label{bothgood}
There are good reasons for studying either diagram (\ref{cd2''}) or diagram (\ref{cd3''}). As we have already seen, one major advantage of diagram (\ref{cd2''}) is that the categories involved are pro-completions. On the other hand, the categories in diagram (\ref{cd3''}) have good exactness properties. Indeed, as $\Pro$ is a regular category, so is $\Pro_{/X}$ by \cite[Example A.5.5]{semiab}, and hence so are $\PGrp_{X}$ and $\PAb_{X}$ by \cite[Theorem 5.11]{barr}. More importantly, $\Pro$ is not exact (see \cite[Example 5.6]{chaus}), but $\PGrp_{X}$ and $\PAb_{X}$ are exact. This essentially follows from \cite[Proposition 5.5.8]{ribesgraph}. In particular, this means that $\PAb_{X}$ is abelian (and in fact satisfies (AB5*) since limits and cokernels are computed on fibres). Actually, by \cite[Remark 2.16(ii)]{jc3}, the category $\PModR_{X}$ is abelian for any fixed profinite ring $R$.
\end{rmk}

\subsection{Main Examples}\label{3.3}

Let us now see various examples, both old and new, where diagram (\ref{cd2''}) or (\ref{cd3''}) appears.

\begin{exmp}\label{maineg}\begin{enumerate}[label=(\roman*)]
\item\label{commoneg1} Consider the original diagram (\ref{cd1}) as written on page \pageref{cd1}. This can be recovered from diagram (\ref{cd3''}) by taking $L'=F\colon\Pro\to\PGrp$ to be the free profinite groups functor. This does arise from Lemma \ref{lem1}, as we saw in Example \ref{exampleJco}\ref{exampleJco1}. If $$(Y,X)=(Y\to X)\in\Pro_{/X},$$ then we conclude from diagram (\ref{cd3''}) that $$F(Y)=\coprod_X{L^*}'(Y)=\coprod_XF(Y(x)),$$ where we recall that $Y(x)$ is the fibre of $Y$ above $x\in X$, and that a topology is implicitly present in the final expression. This is saying that ``the free group over the bundle $Y=\bigsqcup_XY(x)$ is the free product $\coprod_XF(Y(x))$".

If we replace $\PGrp$ with $\PModR$, $\PGrp_{X}=\Grp(\Pro_{/X})$ with $\PModR_{X}$, and the free profinite groups functor with the free profinite $R$-modules functor $R\llbracket-\rrbracket$, then we recover \cite[Proposition A.14]{garethrel}, whose notations are different from ours.

\item\label{commoneg2} Let $k$ be a profinite ring and $H\leq G$ be profinite groups. Consider diagram (\ref{cd2''}) or (\ref{cd3''}) where the right vertical map is the induction functor $$\Ind^G_H(-)\colon\mathbf{PMod}(\kH)\to\mathbf{PMod}(\kG),$$ which is left adjoint to restriction $\Res^G_H(-)$. As we noted at the end of Example \ref{exampleJco}\ref{exampleJco1}, Lemma \ref{lem1} applies here. In particular, we have managed to define induction and restriction in the abelian categories $\mathbf{PMod}(\kH)_{X}$ and $\mathbf{PMod}(\kG)_{X}$ (see Remark \ref{bothgood}). We conclude from diagram (\ref{cd3''}) that given $(P,X)\in\mathbf{PMod}(\kH)_{X}$, we have $$\Ind^G_H\left(\widehat{\bigoplus}_XP(x)\right)=\widehat{\bigoplus}_X(\Ind^G_HP(x)),$$ so ``induction preserves coproducts (internal to $\Pro$)".

More generally, let $R,S$ be profinite rings and $R\to S$ be a map. Then, extension of scalars $(-)\crotimes 
S\colon\PModR\to\mathbf{PMod}(S)$ is left adjoint to restriction $\mathbf{PMod}(S)\to\PModR$, and Lemma \ref{lem1} applies here. By Remark \ref{what}, restriction of scalars preserves $\widehat{\oplus}$ (cf.\ \cite[Proposition 2.8]{jc3}).

\item\label{commoneg3} Let $R,S$ be profinite rings and $M$ be a profinite $R$-$S$-bimodule. Consider diagram (\ref{cd2''}) or (\ref{cd3''}) where the right vertical map is the profinite tensor product $$(-)\crotimes M\colon\mathbf{PMod}(R)\to \mathbf{PMod}(S).$$ This might seem more general than the last part of \ref{commoneg2}, but we cannot in general apply Lemma \ref{lem1} here because the profinite tensor product is not necessarily a left adjoint. However, we can still apply Theorem \ref{mainthm}\ref{mainthm1}. In particular, we have managed to define a tensor product $\mathbf{BPMod}(R)\to \mathbf{BPMod}(S)$. This is a special case of the construction in \cite[Appendix A]{jc3}. To recover the full construction in \cite[Appendix A]{jc3}, we will need to study products of bundle categories, which we will not carry out in the current paper. In addition, restricting to diagram (\ref{cd3''}) defines a tensor product $$(-)\crotimes M\colon\mathbf{PMod}(R)_{X}\to \mathbf{PMod}(S)_{X}$$ where both categories are abelian. We conclude from diagram (\ref{cd3''}) that ``tensor products commute with coproducts".

If $k$ is a profinite ring, $H\leq G$ are profinite groups, $R=\kH$, $S=\kG$, and $M=\kG$, then we leave it to the reader to check that the functor $$\mathbf{PMod}(\kH)_{X}\to\mathbf{PMod}(\kG)_{X}$$ obtained from \ref{commoneg3} is exactly the induction functor obtained from \ref{commoneg2}. Perhaps surprisingly, using the universal property of pro-completions, we have managed to define induction of bundles (which is a left adjoint) \emph{before} defining the tensor product of bundles (which is generally not a left adjoint).

Similarly, for a profinite ring $R$ and a profinite left $R$-module $M$, we can consider diagram (\ref{cd2''}) or (\ref{cd3''}) where the right vertical map is the $\Tor$ functor $\Tor^R_i(-,M)\colon\mathbf{PMod}(R)\to \PAb$. This extends to a $\Tor$ functor on bundles and the commutativity of diagram (\ref{cd3''}) is precisely \cite[Theorem 9.1.1]{ribesgraph}. Some applications of this fact can be found in \cite{jc3}.

\item\label{maineg4} Consider diagram (\ref{cd2''}) or (\ref{cd3''}) where the right vertical map is the abelianisation functor $(-)^\mathrm{ab}\colon\PGrp\to\PAb$, which is left adjoint to the forgetful functor $\PAb\to\PGrp$. Lemma \ref{lem1} trivially applies here, since the abelianisation of a finite group is still finite (unlike in \ref{commoneg1} where a free profinite group/module with finite basis can very much be infinite, and similarly for \ref{commoneg2}). We used here the fact that the functor $G\mapsto G^\mathrm{ab}=G/\overline{[G,G]}$ preserves inverse limits (see \cite[Section IV.2, Exercise 6/Section IV.1, Exercise 5]{ant}), so it is indeed the (unique) extension of the functor $(-)^\mathrm{ab}\colon\Grp_\mathrm{fin}\to\Ab_\mathrm{fin}$. We conclude from diagram (\ref{cd3''}) that given $(P,X)\in\PGrp_{X}$, we have $$\left(\coprod_XP(x)\right)^\mathrm{ab}=\widehat{\bigoplus}_X(P(x)^\mathrm{ab}).$$

To give a concrete example, let $I\in\Set$ and let $X=I\sqcup\{*\}$ be the one-point compactification of $I$. For $i\in I$, let $G_i$ be a profinite group. Then, we can form a bundle $(P,X)\in\PGrp_{X}$, where $P(i)=G_i$ and $P(*)=1$, as follows (see \cite[Example 5.1.1(c)]{ribesgraph}). A local base around $p\in P(i)=G_i$ ($i\neq *$) is given by the open subsets of $G_i$ containing $p$, whereas a local base around $1\in P(*)$ is given by subsets of the form $\bigsqcup_{u\in U}{P(u)}$, where $U$ is an open neighbourhood of $*\in X$, i.e.\ $U$ contains $*$ and almost all elements of $I$. This is called a \emph{bundle of groups over $I$ converging to 1}, and the corresponding coproduct $\coprod_XP(x)\in\PGrp$ is called the \emph{restricted free product of the $G_i$} and denoted by $\coprod^r_{i\in I}G_i$ (see \cite[Example 5.1.9]{ribesgraph}).

We claim that the abelianisation $(P^\mathrm{ab},X)\in\PAb_{X}$ of the above $(P,X)\in\PGrp_{X}$ is also a bundle of (abelian) groups converging to 1. Let us recall how $P^\mathrm{ab}$ is constructed. First write $(P,X)=\varprojlim(P_i,X_i)$, where $(P_i,X_i)\in\BPGrpf$. Then, by definition, $P^\mathrm{ab}\to X$ is the inverse limit of the right maps in the diagrams
\[\begin{tikzcd}
P_i \arrow[rd, two heads] \arrow[rr, two heads] &     & P_i^\mathrm{ab} \arrow[ld, two heads] \\
                                                & X_i &                                      
\end{tikzcd}\]
where $P_i=\bigsqcup_{x\in X_i}P_i(x)$ and $P_i^\mathrm{ab}=\bigsqcup_{x\in X_i}(P_i(x)^\mathrm{ab})$. The top map in the inverse limit
\[\begin{tikzcd}
P \arrow[rd, two heads] \arrow[rr, two heads] &   & P^\mathrm{ab} \arrow[ld, two heads] \\
                                              & X &                                    
\end{tikzcd}\]
when restricted to each fibre is just the usual abelianisation $P(x)\sur P^\mathrm{ab}(x)=P(x)^\mathrm{ab}$. In particular, the top map is a quotient map (which also follows from the simple fact that $P\sur P^\mathrm{ab}$ is surjective). We leave it to the reader to verify that $(P^\mathrm{ab},X)$ is therefore also a bundle converging to 1.

By \cite[Example 5.6.4(c)]{ribesgraph}, the internal coproduct $\widehat{\bigoplus}_X(P(x)^\mathrm{ab})$ of $(P^\mathrm{ab},X)\in\PAb_{X}$ is the categorical \emph{product} $\prod_{i\in I}P(i)^\mathrm{ab}=\prod_{i\in I}G_i^\mathrm{ab}\in\PAb$. We therefore obtain the unintuitive statement $$\left(\coprod^r_{i\in I}G_i\right)^\mathrm{ab}=\prod_{i\in I}G_i^\mathrm{ab}.$$ However, as a special case, we recover the more reasonable statement that the abelianisation of $$F(X)=F(I\sqcup\{*\})=\coprod^r_I\widehat{\Z},$$ the free profinite group on $I$ (converging to 1), is $$F(I\sqcup\{*\})^\mathrm{ab}=\prod_I\widehat{\Z},$$ the free proabelian group on $I$ (converging to 1): see \cite[Example 5.1.3(c)]{ribesgraph} and \cite[Example 3.3.8(c)]{profinite}.
\end{enumerate}\end{exmp}

\appendix
\section{Pontryagin duality for bundles}

This paper is perhaps incomplete without at least mentioning Pontryagin duality, which forms a significant part of \cite{gareth} and \cite{boggi} (see \cite[Theorem 4.1]{gareth} and \cite[Section 2.6]{boggi}). The paper \cite{gareth} utilises sheaf theory to obtain Pontryagin duality for bundles, but we will directly describe Pontryagin duality in terms of bundles. The following is therefore not original and can be recovered by carefully studying \cite{gareth}.

Let us first point out that the perspectives of \cite{haran} and \cite{melfree} (and therefore also of \cite{ribesgraph}, \cite{gareth} and ours) are equivalent, as mentioned in the addendum of \cite{melfree}. On the other hand, the categories under consideration in \cite{boggi} are slightly more restrictive, as explained to the author by Boggi.

Let $R$ be a fixed profinite ring and consider the category $\BPModRf$ of bundles of finite (discrete topological) $R$-modules and its pro-completion $\BPModR$. Recall that an object of $\BPModRf$ is a pair $(P,X)$ of finite sets together with a function $p\colon P\to X$, such that each fibre $P(x)=p^{-1}(\{x\})$ is a finite $R$-module. A morphism from $(P,X)$ to $(Q,Y)$ consists of maps $a\colon X\to Y$ and $f\colon P\to Q$ making the obvious square commute, such that the restriction of $f$ to each $x\in X$ is a module homomorphism $P(x)\to Q(a(x))$.

The dual category will be given by the category $\mathbf{B}^\mathrm{op}\ModR_\mathrm{fin}$ of ``opposite bundles of finite $R$-modules", which one obtains by making all the spaces finite in \cite[Definition 2.1]{gareth}. Spelt out, $\mathbf{B}^\mathrm{op}\ModR_\mathrm{fin}$ has the same objects as $\BPModRf$, but a morphism from $(P,X)$ to $(Q,Y)$ now consists of maps $b\colon Y\to X$ and $g\colon P\times_X Y\to Q$ making the triangle
\[\begin{tikzcd}
P\times_X Y \arrow[rd, "\pi_Y"'] \arrow[rr, "g"] &   & Q \arrow[ld] \\
                                                 & Y &             
\end{tikzcd}\]
commute, such that the restriction of $g$ to each $y\in Y$ is a module homomorphism $P(b(y))\to Q(y)$.

There is a functor $$(-)^\vee=\Hom(-,\mathbb{Q}/\mathbb{Z})\colon\BPModRf^\mathrm{op}\to\mathbf{B}^\mathrm{op}\mathbf{Mod}(R^\mathrm{op})_\mathrm{fin},$$ the \emph{Pontryagin duality functor}, which sends $(P,X)\in\BPModRf^\mathrm{op}$ to the ``opposite bundle" $(P^\vee,X)$, where $P^\vee(x)=P(x)^\vee$. Given a morphism $(P,X)\to(Q,Y)$ in $\BPModRf^\mathrm{op}$ (i.e.\ maps $b\colon Y\to X$ and $f\colon Q\to P$ making the obvious square commute, such that each $Q(y)\to P(b(y))$ is a module homomorphism), the functor $(-)^\vee$ sends it to the morphism $(P^\vee,X)\to(Q^\vee,Y)$ in $\mathbf{B}^\mathrm{op}\mathbf{Mod}(R^\mathrm{op})_\mathrm{fin}$, where $b\colon Y\to X$ is as above, and $g\colon P^\vee\times_XY\to Q^\vee$ is given on $y\in Y$ by taking the usual Pontryagin dual of $Q(y)\to P(b(y))$.

Similarly, there is a functor $(-)^\vee\colon\mathbf{B}^\mathrm{op}\mathbf{Mod}(R^\mathrm{op})_\mathrm{fin}\to\BPModRf^\mathrm{op}$. The following is then easy to verify.

\begin{lemma}
The functor $(-)^\vee=\Hom(-,\mathbb{Q}/\mathbb{Z})\colon\BPModRf^\mathrm{op}\to\mathbf{B}^\mathrm{op}\mathbf{Mod}(R^\mathrm{op})_\mathrm{fin}$ is an equivalence, i.e.\ it establishes a duality between $\BPModRf$ and $\mathbf{B}^\mathrm{op}\mathbf{Mod}(R^\mathrm{op})_\mathrm{fin}$.
\end{lemma}

Taking pro-completions on both sides, we then get the equivalence $$\BPModR=\Pro(\BPModRf)=\Pro(\mathbf{B}^\mathrm{op}\mathbf{Mod}(R^\mathrm{op})^\mathrm{op}_\mathrm{fin})=\mathbf{Ind}(\mathbf{B}^\mathrm{op}\mathbf{Mod}(R^\mathrm{op})_\mathrm{fin})^\mathrm{op},$$ where $\mathbf{Ind}$ means ind-completion. This therefore establishes a (Pontryagin) duality between $\BPModR=\Pro(\BPModRf)$ and $\mathbf{Ind}(\mathbf{B}^\mathrm{op}\mathbf{Mod}(R^\mathrm{op})_\mathrm{fin})$.

Next, let us prove that the Pontryagin duality functor we have defined is in fact equivalent (i.e.\ naturally isomorphic) to the one in \cite{gareth}. Since we have extended $(-)^\vee$ using the universal property of pro-completions, it suffices to show that the corresponding functor in \cite{gareth} commutes with inverse limits and agrees with ours on finite bundles. But the functor in \cite{gareth} commutes with inverse limits (i.e.\ sends inverse limits in $\BPModR$ to direct limits in $\mathbf{Ind}(\mathbf{B}^\mathrm{op}\mathbf{Mod}(R^\mathrm{op})_\mathrm{fin})$), by the Pontryagin duality established in \cite{gareth}. For finite bundles, it is easy to check directly that the Pontryagin duality functor in \cite{gareth} agrees with ours. Indeed, we only need to know that the Pontryagin duality functor in \cite{gareth} commutes with taking fibres: see \cite[Theorem 4.3]{gareth}.

We should point out there is a dual theory to what we have done for modules in Section \ref{sec3}, which is studied in \cite{gareth}, \cite{boggi} and \cite{jc3}. We leave it to the interested reader to figure out the details.

\bibliographystyle{unsrt}

\begin{thebibliography}{10}

\bibitem{haran}
Dan Haran.
\newblock On closed subgroups of free products of profinite groups.
\newblock {\em Proceedings of the London Mathematical Society}, 3(2):266--298,
  1987.

\bibitem{melfree}
Oleg~Vladimirovich Mel'nikov.
\newblock Subgroups and the homology of free products of profinite groups.
\newblock {\em Izv. Akad. Nauk SSSR Ser. Mat.}, 53(1):97--120, 1989.

\bibitem{ribesgraph}
Luis Ribes.
\newblock {\em Profinite graphs and groups}, volume~66.
\newblock Springer, 2017.

\bibitem{gareth}
Gareth Wilkes.
\newblock Pontryagin duality and sheaves of profinite modules.
\newblock {\em Journal of Algebra}, 2025.

\bibitem{boggi}
Marco Boggi.
\newblock Lannes' ${T}$-functor and mod-$ p $ cohomology of profinite groups.
\newblock {\em Journal of Algebra}, 2026.

\bibitem{jc3}
Jiacheng Tang.
\newblock Profinite direct sums with applications to profinite groups of type
  ${\Phi}_{R}$.
\newblock {\em Bulletin of the London Mathematical Society}, 2025.

\bibitem{garethrel}
Gareth Wilkes.
\newblock Relative cohomology theory for profinite groups.
\newblock {\em J. Pure Appl. Algebra}, 223(4):1617--1688, 2019.

\bibitem{maccat}
Saunders Mac~Lane.
\newblock {\em Categories for the working mathematician}, volume~5.
\newblock Springer Science \& Business Media, 2013.

\bibitem{profinite}
Luis Ribes and Pavel Zalesskii.
\newblock {\em Profinite groups}.
\newblock Springer, 2000.

\bibitem{handbook}
Francis Borceux.
\newblock {\em Handbook of categorical algebra: Basic category theory},
  volume~1.
\newblock Cambridge University Press, 1994.

\bibitem{elephant}
Peter~T Johnstone.
\newblock {\em Sketches of an Elephant: A Topos Theory Compendium}, volume~2.
\newblock Oxford University Press, 2002.

\bibitem{catlogic}
Bart Jacobs.
\newblock {\em Categorical logic and type theory}, volume 141.
\newblock Elsevier, 1999.

\bibitem{ribesamal}
Luis Ribes.
\newblock On amalgamated products of profinite groups.
\newblock {\em Mathematische Zeitschrift}, 123(4):357--364, 1971.

\bibitem{venj}
Otmar Venjakob.
\newblock Subgroup theorems for free profinite products with amalgamation.
\newblock {\em Journal of the London Mathematical Society}, 68(1):12--24, 2003.

\bibitem{trees}
Jean-Pierre Serre.
\newblock {\em Trees}.
\newblock Springer Science \& Business Media, 2002.

\bibitem{gro1}
Michael Artin, Alexander Grothendieck, and Jean-Louis Verdier.
\newblock Th{\'e}orie des topos et cohomologie {\'e}tale des sch{\'e}mas.
\newblock {\em Lecture notes in mathematics}, 1972.

\bibitem{ulmer}
Friedrich Ulmer.
\newblock Properties of dense and relative adjoint functors.
\newblock {\em Journal of Algebra}, 8(1):77--95, 1968.

\bibitem{lawvere}
F~William Lawvere.
\newblock Functorial semantics of algebraic theories.
\newblock {\em Proceedings of the National Academy of Sciences},
  50(5):869--872, 1963.

\bibitem{proadj}
{nLab authors}.
\newblock ``pro-left adjoint", 2021.
\newblock URI: https://ncatlab.org/nlab/show/pro-left+adjoint.

\bibitem{semiab}
Francis Borceux and Dominique Bourn.
\newblock {\em Mal'cev, protomodular, homological and semi-abelian categories},
  volume 566.
\newblock Springer Science \& Business Media, 2004.

\bibitem{barr}
Michael Barr.
\newblock {\em Exact categories}.
\newblock Springer, 1971.

\bibitem{chaus}
Vincenzo Marra and Luca Reggio.
\newblock A characterisation of the category of compact {H}ausdorff spaces.
\newblock {\em Theory Appl. Categ.}, 35:Paper No. 51, 1871--1906, 2020.

\bibitem{ant}
J{\"u}rgen Neukirch.
\newblock {\em Algebraic number theory}, volume 322.
\newblock Springer Science \& Business Media, 2013.

\end{thebibliography}

\end{document}